\documentclass[12pt]{article}
\usepackage{fancyhdr, array,calc,graphicx,url,tabularx}
\usepackage{geometry}
\usepackage{latexsym}
\usepackage{amssymb}
\usepackage{amsmath}
\usepackage{makeidx}

\usepackage{enumitem,ifthen}
\usepackage{xcolor,xspace}
\usepackage{tikz}
\usetikzlibrary{matrix,arrows,cd,positioning,scopes}%

\tikzset{
	  symbol/.style={
		      draw=none,
		          every to/.append style={
				        edge node={node [sloped, allow upside down, auto=false]{$#1$}}}
				  }
				  }
				  \usepackage[hidelinks,colorlinks]{hyperref}

				  \usepackage{graphicx} % a package for including diagrams
				  \usepackage{setspace} % provides the doublespace command below
				  \usepackage{latexsym} % must have basic LaTeX symbols
				  \usepackage{amsfonts} % another must have AMS fonts
				  \usepackage{amsmath} % AMS math stuff like matrices and such
				  \usepackage{amssymb}
				  \usepackage{accents}
				  \usepackage{bm}
				  \usepackage{graphics}

				  \makeindex

\newcommand{\cF}{{\mathcal F}}

\newcommand{\cM}{{\mathcal M}}

\newcommand{\cT}{{\mathcal T}}
\newcommand{\cU}{{\mathcal U}}

\newcommand{\R}{{\mathbb R}}

\newcommand{\adp}{\mbox{{\sf AD}}^+}

\newcommand{\ad}{\mbox{{\sf AD}}}

\newcommand{\adr}{\mbox{{\sf AD}}_{\mathbb R}}

\newcommand{\hd}{\mbox{HOD}}

\newcommand{\hpc}{\mbox{{\sf HPC}}}

\newcommand{\meas}{\mbox{meas}}

\newcommand{\hc}{\mbox{HC}}

\newcommand{\dom}{\mathop{\mathrm{dom}}\nolimits}

\newcommand{\cf}{\mathop{\mathrm{cof}}\nolimits}

\newcommand{\lh}{\mathop{\mathrm{lh}}\nolimits}

\newcommand{\crit}{\mathop{\mathrm{crit}}\nolimits}
\newcommand{\ran}{\mathop{\mathrm{ran}}\nolimits}

\newcommand{\ult}{\mathop{\mathrm{Ult}}\nolimits}

\newcommand{\cd}{\mathop{\mathrm{Code}}\nolimits}

\usepackage{amsthm} %formatting theorems etc

\newtheorem{theorem}{Theorem}[section]
\newtheorem{lemma}[theorem]{Lemma}
\newtheorem{corollary}[theorem]{Corollary}

\theoremstyle{definition}
\newtheorem{definition}[theorem]{Definition}

\newtheorem{claim-plain}[theorem]{Claim}

\theoremstyle{remark}

\newcommand{\conc}{{}^\smallfrown} %concatenation symbol

\title{Suslin cardinals and cutpoints in mouse limits}
\author{Stephen Jackson, Grigor Sargsyan, and John Steel }
\date{July 2022}

\begin{document}

\maketitle

\section{Introduction}

     We assume $\adp$ throughout, and assume familiarity with the main definitions and
     results of \cite{nitcis} and \cite{mouse.suslin} concerning mouse pairs $(P,\Sigma)$ and their
     associated mouse limits $M_\infty(P,\Sigma)$. By way of a brief review:
     a {\em mouse pair} consists of
     a countable pfs premouse $P$ together with an iteration strategy $\Sigma$ for $P$
     having certain regularity properties. Here ``pfs" stands for "projectum-free spaces",
    and corresponds to a minor variation on the usual Jensen-indexed fine structure.
    Our pfs premice will always be {\em projectum stable}, where
    $P$ is projectum stable iff $P$ has type 1\footnote{Cf. \cite[\S 4.1]{nitcis}},
and for $k=k(P)$ the distinguished soundness degree of $P$,
the $r\Sigma_k$ cofinality of $\rho_k(P)$ is not measurable in $P$.\footnote{ This property
is called {\em strong stability} in \cite{nitcis}. It holds trivially if $P$ has type 1 and
$k(P)=0$.}
$M_\infty(P,\Sigma)$ is the direct limit
of all countable $\Sigma$-iterates of $P$. If $P$ is projectum stable, then
$M_\infty(P,\Sigma)$ exists and is itself projectum stable.
     
     Our motivation is the following conjecture.

\bigskip
\noindent
{\em Conjecture.} Let $(P,\Sigma)$ be a projectum stable mouse pair, and let
$\kappa$ be a cardinal of $V$ such that $\kappa < o(M_\infty(P,\Sigma))$;
then the following are equivalent:
\begin{itemize}
	\item[(1)] $\kappa$ is a Suslin cardinal,
	\item[(2)] $\kappa$ is a cutpoint of $M_\infty(P,\Sigma)$.
\end{itemize}

The conjecture would imply that assuming $\hpc$\footnote{$\hpc$ is the assertion that
the sets of reals coding least branch $\hd$ pairs are Wadge-cofinal in the Suslin-co-Suslin
sets of reals.}, the Suslin cardinals
are precisely the cardinals of $V$ that are cutpoints of the extender sequence of $\hd$.
That (2) implies (1) follows easily from work in \cite{mouse.suslin}:

\begin{lemma}\label{twoimpliesone} Let $(P,\Sigma)$ be a projectum stable mouse pair,
	and $\kappa < o(M_\infty(P,\Sigma))$ be a cardinal of $V$.
	Suppose $\kappa$ is a cutpoint of $M_\infty(P,\Sigma)$; then
	$\kappa$ is a Suslin cardinal.
\end{lemma}
We shall prove the lemma below.

Recent work of Jackson and Sargsyan gets us a lot closer to a proof of 
the converse direction. The part of this paper that goes beyond
\cite{mouse.suslin} is mostly an account of their work.

\begin{definition} For any $\kappa$, $\meas_{\kappa}$ is the collection
of all ultrafilters on $\kappa$
\end{definition}
Of course, each $U \in \meas_\kappa$ is countably complete, and Rudin-Kielser
reducible to the Martin measure on degrees if $\kappa < \theta$.

G. Sargsyan recently proved the following.

\begin{theorem}[Sargsyan \cite{paulgrigor}]\label{paulgrigor}
Assume $\adp$, and let $(P,\Sigma)$ be a projectum stable mouse pair. Let 
$E$ be an extender of the sequence of $M_\infty(P,\Sigma)$ with critical
point $\kappa$, and 
such that 
\begin{itemize}
\item[(1)] $\kappa$ is a cutpoint of $M_\infty(P,\Sigma)$, and
\item[(2)] $E$ is total on $M_{\infty}(P,\Sigma)$, and
	$\kappa < \rho_n(M_\infty(P,\Sigma))$, for $n = k(P)$.
\end{itemize}
	Then there is a  $U \in \meas_{\kappa}$ 
such that if $j_U \colon V \to \ult(V,U)$ 
and $i_E \colon M_\infty(P,\Sigma) \to
\ult(M_\infty(P,\Sigma),E)$ are the canonical embeddings, then
\[
j_U \restriction M_\infty(P,\Sigma) = \sigma \circ i_E,
\]
for some elementary $\sigma \colon \ult(M_\infty(P,\Sigma),E) \to
j_U(M_\infty(P,\Sigma))$, and hence
	\[
		\lambda_E \le j_U(\kappa).
		\]
\end{theorem}
See \cite{paulgrigor}. (\cite{paulgrigor} does not state the result
in this generality, but this is what the proof gives.)  Using the
known connections between Suslin cardinals, measures, and Martin classes
(see \cite{jackson}[\S 3]), Theorem \ref{paulgrigor} yields at once
Steel's theorem that under $\adr + \hpc$, every point $\theta_\alpha$
in the Solovay sequence is a cutpoint of the $\hd$-sequence.
\footnote{See \cite[Theorem 5.1]{mouse.suslin} It is important here that we are talking about
cutpoints with respect to extenders on the $\hd$-sequence. We do not 
have a proof that every $\theta_{\alpha}$ is a cutpoint with respect to extenders
belonging to $\hd$. The results of \cite{farmer} would seem to be
relevant there.} The resulting proof is simpler and
more general than that of \cite{mouse.suslin}.

Jackson has recently observed
that the results on Martin classes of \cite{jackson}[\S 3] can be extended so as
to prove the following. 

\begin{theorem}[Jackson]\label{measurebound} Assume $\adp$.
Let $\kappa$ be a limit of Suslin cardinals of uncountable cofinality, and
$\lambda$ the least Suslin cardinal $>\kappa$; then for any ultrafilter
$U$ on $\kappa$, 
$j_U(\kappa) < \lambda$.
\end{theorem}

We remark that if $\lambda$ is Suslin, then $\lambda \le
\sup(\lbrace j_U(\kappa) \mid \kappa < \lambda \wedge U \in \meas_{\kappa} \rbrace)$.
This comes from Martin-Solovay construction of a scale on $\neg p[T]$, where
$T$ is weakly homogeneous. We do not know whether the reverse inequality holds
at all Suslin cardinals $\lambda$. If $\lambda$ is a limit of Suslins, it is
trivial. If $\lambda$ is the next Suslin after a limit of Suslins, then the theorem above
says a lot, but does not fully answer the question.

\begin{definition} \label{oq} For any premouse $Q$ and $\kappa <o(Q)$,
	\begin{itemize}
		\item[(a)] 
$o(\kappa)^Q$ is the strict sup of all $\eta$ such that
$\crit(E_\eta^Q)=\kappa$. If there are no such $\eta$,
			then $o(\kappa)^Q = 0$.
			\item[(b)] $\kappa$ is {\em $Q$-regular} iff
				there is no $\eta < \kappa$ and total $\Sigma_{k(Q)}(Q)$
				function $f \colon \eta \to \kappa$ with range
				cofinal in $\kappa$.
			\item[(c)] $\kappa$ is {\em $Q$-measurable} iff
				$\kappa < \rho_{k(Q)}(Q)$, and 
				$o(\kappa)^Q \eta \ge (\kappa^+)^Q$.
	\end{itemize}
\end{definition}
Coherence implies that if
$\kappa$ is a cutpoint of $Q$, so is $o(\kappa)^Q$. Our notion of regularity
involves all functions that might be used in some nondropping ultrapower
of $Q$. Thus we have by the usual ``regulars are measurable" argument 

\begin{lemma}\label{regularsmeasurable} Let $(P,\Sigma)$ be a projectum stable mouse pair,
$M_\infty = M_\infty(P,\Sigma)$, and $\kappa < o(M_\infty)$ have
uncountable cofinality in $V$; then $\kappa$ is $M_\infty$-regular
iff $\kappa$ is $M_\infty$-measurable.
\end{lemma}

	Putting the two theorems above together, with some sauce from \cite{mouse.suslin},
we get the following.

\begin{theorem} \label{addsauce}  Assume $\adp$, let $(P,\Sigma)$
be a projectum stable mouse pair, and let $M_\infty = M_\infty(P,\Sigma)$.
Let $\kappa < o(M_\infty)$ be a limit of Suslin
cardinals such that
$\cf(\kappa)>\omega$ in $V$; then
	\begin{itemize}
		\item[(1)]$\kappa$ is a limit of cutpoints of $M_\infty$, and
	\end{itemize} Suppose $\kappa^+ \le o(M_\infty)$, and let 
$\lambda$ be the least Suslin cardinal $>\kappa$; then
\begin{itemize}
	\item[(2)] there is a cutpoint $\mu$ of $M_\infty$ such that
		$\kappa \le \mu < (\kappa^+)^V$ and $o(\mu)^{M_\infty} = \lambda$,
	\item[(3)] if $o(\kappa)^{M_\infty} \ge (\kappa^+)^V$, then $o(\kappa)^{M_\infty} = \lambda$, and 
	\item[(4)] if $S_\kappa$ is closed under $\forall^{\R}$, then $o(\kappa)^{M_\infty} = \lambda$.
\end{itemize}
\end{theorem}
Some of these results were proved in \cite{mouse.suslin} in the case that $(P,\Sigma)$ is a
pointclass generator.

From this we get at once

\begin{corollary} \label{hodcorollary} Assume $\adr + \hpc$, and let $\kappa$ be a limit
of Suslin cardinals of uncountable cofinality, and regular in $\hd$,
and let $\lambda$ be the least
Suslin cardinal $>\kappa$; 
then
\begin{itemize}
	\item[(1)] $(\kappa^+)^{\hd} \le o(\kappa)^{\hd} \le \lambda$,
	\item[(2)] there is a cutpoint $\mu$ of $\hd$ such that
$\kappa \le \mu < (\kappa^+)^V$ and $o(\mu)^{\hd} = \lambda$,
\item[(3)] if $o(\kappa)^{\hd} \ge (\kappa^+)^V$, then $o(\kappa)^{\hd} = \lambda$, and
	        \item[(4)] if $S_\kappa$ is closed under $\forall^{\R}$, then $o(\kappa)^{\hd} = \lambda$.
\end{itemize} 
\end{corollary}

If $\kappa$ is a countable cofinality limit of Suslins, then $\kappa^+$
is the next Suslin (and somewhat like $\omega_1$). See \cite{jackson}[3.28].
In this case we have

\begin{theorem}\label{cofomegalimit} 
Let $(P,\Sigma)$ be a projectum stable mouse pair, $M_\infty =
M_\infty(P,\Sigma)$, and $\kappa$ be a limit of
Suslin cardinals of countable $V$-cofinality.
Suppose
$(\kappa^+)^V \le o(M_\infty)$; then 
	\begin{itemize}
		\item[(1)] $\kappa$ and $(\kappa^+)^V$ are limits of cutpoints
			in $M_\infty$, and
		\item[(2)] $(\kappa^+)^V < \rho_{k(P)}(M_\infty)$, and
	$(\kappa^+)^V$ is the critical point of a total extender
			from the $M_\infty$-sequence.
	\end{itemize}
\end{theorem}
And then of course there is a corollary for $\hd$  parallel to \ref{hodcorollary}.

These results seem close to a proof of the conjecture. What's missing is the analog of Jackson's
result on measure-bounding for the Suslin
cardinals corresponding to higher levels of a projective-like hierarchy. For the
ordinary projective hierarchy, Jackson has proved these as part of his computation
of the projective ordinals. But perhaps the full force of this machinery is not needed.
We do have

\begin{theorem}\label{conjectureprojective} The conjecture holds when $\kappa$ is one of
	the $\bf{\delta}^1_{2n+1}$'s or their cardinal predecessors.
\end{theorem}

Theorem \ref{conjectureprojective} was known for various natural $(P,\Sigma)$
by other means already. Theorem \ref{paulgrigor}, \cite{mouse.suslin}, and
Jackson's results on measure
bounding in the projective hierarchy yield a different, more general proof.\footnote{So for example,
when $P$ is $M_1$ cut at its Woodin, and $\Sigma$ is its canonical strategy,
we get a new proof that $o(M_\infty(P,\Sigma)) = \aleph_\omega$, with the least strong of $M_\infty$ being
$< \omega_2$. And in fact, something similar must happen for any $(P,\Sigma)$ such that
$o(M_\infty(P,\Sigma)) \ge \omega_2$.}

In this note we shall prove the results above.

\section{ Proof of Theorem \ref{paulgrigor}}

 Let $(P,\Sigma)$ be a projectum stable mouse pair, $\cF(P,\Sigma)$ the directed system of all its nondropping
 iterates, and  $M_\infty = M_\infty(P,\Sigma)$ the direct limit of $\cF(P,\Sigma)$.
 For the associated iteration maps of the system we write $\pi_{Q,R} \colon Q \to \R$ and
 $\pi_{Q,\infty} \colon Q \to M_\infty$. It's ok here to drop mention of the 
 strategy of $Q$, since we are dealing exclusively with tails of a single positional strategy $\Sigma$.
 Let $k = k(P)$. \footnote{This is
 the quantifier level at which we are condisdering $P$. $P$ is always $k(P)$ sound, but it may not be
 $k(P)+1$ sound. See \cite{nitcis}[Chapter 2].} Let $E$
 be an extender on the sequence of $M_\infty$, and $\kappa = \crit(E)$. Suppose that $\kappa$ is a cutpoint
 of $P$ (and hence a limit of cutpoints), and that $E$ is total on $P$ and $\kappa < \rho_k(P)$.
 We want to embed $\ult(M_\infty,E)$ into $j_U(M_\infty)$, for some ultrafilter $U$ on $\kappa$.

By replacing $(P,\Sigma)$ with an iterate of itself, we may assume $E \in \ran(\pi_{P,\infty})$.
For any $(R,\Sigma_R) \in \mathcal{F}(P,\Sigma)$, let
\[
	\pi_{R,\infty}(E_R) = E,
	\]
	and
	\[
		\pi_{R,\infty}(\kappa_R) = \kappa.
		\]

 If $d$ is
 a Turing degree and $Q \in \hc$, we write $Q \le d$ to mean that $Q$ is coded by a real recursive in $d$.
 (Fix some natural coding system.) Let
 \[
	 \cF(d) = \lbrace (Q,\Sigma_Q) \mid Q \le d \wedge (Q,\Sigma_Q) \in \cF(P,\Sigma)\rbrace,
	 \]
	 and
	 \[
		 M_d = \text{ result of simulaneously comparing all $(Q,\Sigma_Q) \in \cF(d)$}.
		 \]
We note here that by \cite{siskindsteel}, $\Sigma$ is positional. It follows that comparisons
between iterates of $(P,\Sigma)$ never encounter strategy disagreements, and so can be done
by iterating away least extender disagreements as usual. The simultaneous comparsion
referred to above proceeds by iterating away least extender disagreements. It does not
depend on any enumeration of $d$, just $d$ itself. Set also
\[
	\Sigma_d = \Sigma_{M_d},
	\]
so that $(M_d,\Omega_d) \in \cF(P,\Sigma)$.  
Let $(E_d, \kappa_d) = (E_{M_d},\kappa_{M_d})$, and
\[
Q_d = M_d||\lh(E_d),
\]
and
\[
\Omega_d = (\Sigma_d)_{Q_d}
\]
be the strategy for $Q_d$ that is part of $\Sigma_d$. Our ``$||$" notation
indicates that $Q_d$ is passive, that is, the last extender predicate $E_d$
has been removed. 

\medskip
\noindent
{\em Claim 1.} $M_\infty(Q_d,\Omega_d) \unlhd M_\infty$.

\medskip
\noindent
{\em Proof.} Let $R = Ult(M_d,E_d)$ and $S= \ult(R,F)$, where
$F$ is the order zero total measure on $\lambda(E_d)=i_{E_d}^{M_d}(\kappa_d)$.
Then $(S,\Sigma_S) \in \cF(P,\Sigma)$, and $(Q_d,\Omega_d)$ is a cardinal
cutpoint initial segment of $(S,\Sigma_S)$. The claim follows.
\hfill    $\square$

\medskip

We can now define our ultrafilter $U$ on $\kappa$. Let
$\bar{\pi}_{d,\infty} = \pi_{Q_d,\infty}^{\Omega_d}$. For $A \subseteq \kappa$,
\[
	A \in U \text{ iff } \forall^*d (\bar{\pi}_{d,\infty}(\kappa_d) \in A).
	\]
Here $\forall^*d$ refers to the Martin measure. $U$ is clearly a countable complete ultrafilter on $\kappa$.
We must now define the desired $\sigma \colon \ult(M_\infty,E) \to j_U(M_\infty)$.
Of course, the definition will be in the form
$\sigma(i_E^{M_\infty}(f)(a)) = j_U(f)(\sigma(a))$. We just have
to figure out what $\sigma(a)$ is. 

Fix an $a \in [\lambda(E)]^{<\omega}$, and let $(R,\Sigma_R)$
be such that $a \in \ran(\pi_{R,\infty})$. Say
\[
	\pi_{R,\infty}(a_R) =a.
	\]
We shall define a function
$f_R^a$, and show that $[f_R^a]_U$ is independent of the $R$ we have
chosen. We then set $\sigma(a) = [f_R^a]_U$. Towards defining
$f^a_R$, let $d$ be any degree such that $R \le d$,
and set
\[
	a_d = \pi_{R,M_d}(a_R) =\pi_{M_d,\infty}^{-1}(a).
	\]
	
The main claim is the following. 

\bigskip
\noindent
{\em Claim 2.} Let $a \in [\lambda(E)]^{<\omega}$, and suppose
$ R \le c$ and $R \le d$.
Suppose $\bar{\pi}_{c,\infty}(\kappa_c) = \bar{\pi}_{d,\infty}(\kappa_d)$;
then $\bar{\pi}_{c,\infty}(a_c) = \bar{\pi}_{d,\infty}(a_d)$.

\medskip
\noindent
{\em Proof.} 
Let $\cT$ be the normal tree by $\Sigma_R$ from $R$ to $M_c$,
and let $\alpha$ be least
such that $\lh(E_\alpha^{\cT}) \ge \lh(E_c)$. Since $E_c$ is on the
last model of $\cT$, $\lh(E_\alpha^{\cT})> \lh(E_c)$. Also,
$E_c \in \ran(i_{0,\infty}^{\cT})$, so $\alpha$ is on the main
branch of $\cT$, and $\crit(i_{\alpha,\infty}^{\cT}) > \lh(E_c)$.
Note that $Q_c \unlhd \cM_\alpha^\cT$.

Similarly, let $\cU$ be the normal tree by $\Sigma_R$ from $R$ to $M_d$,
and let $\beta$ be least
such that $\lh(E_\beta^{\cU}) \ge \lh(E_d)$. Since $E_d$ is on the
last model of $\cU$, $\lh(E_\beta^{\cU})> \lh(E_d)$. Again, 
$\beta$ is on the main
branch of $\cU$, $\crit(i_{\beta,\infty}^{\cU}) > \lh(E_d)$,
and
$Q_d \unlhd \cM_\beta^\cU$.

Now notice that $M_\infty(Q_c,\Omega_c) = M_\infty(Q_d,\Omega_d)$. 
This is because both are cutpoint initial segments of $M_\infty(P,\Sigma)$,
and both have a top block that begins at the same place, namely
$\bar{\pi}_{c,\infty}(\kappa_c)
= \bar{\pi}_{d,\infty}(\kappa_d)$. It follows that $(Q_c,\Omega_c)$ is mouse equivalent
to $(Q_d,\Omega_d)$ (see \cite{mouse.suslin}[2.2]). They compare by iterating
away least extender disagreements, because we are working with tails of a single
positional strategy. Let $\cT_1$ on $Q_c$ and $\cU_1$ on $Q_d$ be the normal trees
with common last model $S$ that we get from this comparison. It is enough to see
that
\[
	i^{\cT_1}(a_c) = i^{\cU_1}(a_c),
	\]
where these are the main branch embeddings of $\cT_1$ and $\cU_1$. (Note
$i^{\cT_1} = \pi_{Q_c,S}^{\Omega_c}$ and $i^{\cU_1} = \pi_{Q_d,S}^{\Omega_d}$.)

To see this, consider the normal trees 
\[
	\cT_0 = \cT \restriction (\alpha+1) \conc \langle E_c \rangle
	\]
	and
	\[
		\cU_0 = \cU \restriction (\beta+1) \conc \langle E_d \rangle.
	\]
	It is important here that we are talking about normal extensions; $E_c$
may not be applied to $\cM_\alpha^{\cT}$, but instead some earlier model.
Letting $N_0 = \cM_{\alpha+1}^{\cT_0}$ and $N_1 = \cM_{\beta+1}^{\cU_0}$ be the last
models, we have that $Q_c$ is a cardinal cutpoint initial segment of $N_0$, and
$o(Q_c) < \rho_k(N_0)$, and similarly for $Q_d$ and $N_1$. Thus $\cT_1$ and $\cU_1$
can be considered as normal, nondropping\footnote{The reader must have figured out by now
that this means the main branch does not drop.} trees on $N_0$ and $N_1$. Let us do that. Let
\[
	X = X(\cT_0,\cT_1),
\]
and
\[
	Y = X(\cU_0,\cU_1)
\]
be the full normalizations of the two stacks, so that $X$ and
$Y$ are normal trees on $R$ by $\Sigma_R$. (See \cite{localhodcomp} or
\cite{siskindsteel}.) We can write
\[
	X = X_0 \conc \langle F \rangle
	\]
	and
	\[
		Y = Y_0 \conc \langle G \rangle,
\]
where $X_0$ and $Y_0$ have last models $N_0^*$ and $N_1^*$ respectively,
both extend $S$, and $F$ and $G$ are the extenders with index $o(S)$
in the two models. Now note that the generators for the branch extender
$R$-to-$N_0^*$ in $X_0$ are contained in $o(S)$, as are the generators
of $R$-to-$N_1^*$ in $Y_0$. So both are trees by $\Sigma_R$ using only
extenders of length $<o(S)$, so in fact,
\[
	X_0 = Y_0
\]
and $N_0^* = N_1^*$, and $i^{X_0} = i^{Y_0}$. 

Let $\Phi$ and $\Psi$ be the weak tree embeddings of
$\cT_0$ and $\cU_0$ into $X$ and $Y$ that come from
full normalization. We have
\[
	t_{\alpha +1}^\Phi \colon N_0 \to N_0^*
	\]
	and
	\[
		t^\Psi_{\beta+1} \colon N_1 \to N_1^*,
		\]
from that process, with
\[
	t_{\alpha+1}^\Phi \restriction \lh(E_c) = t_\alpha^\Phi \restriction \lh(E_c) =
i^{\cT_1} \restriction \lh(E_c),
\]
and
\[
t_{\beta+1}^\Psi \restriction
\lh(E_d) = t_\beta^\Psi = i^{\cU_1} \restriction \lh(E_d).
\]
Also,
\[
	i^{X_0} = t_{\alpha}^\Phi \circ i^{\cT_0}_{0,\alpha}
\]
and 
\[
	i^{Y_0} = t_{\beta}^\Psi \circ i^{\cU_0}_{0,\beta},
	\]
by the way normalization works. Letting $a_R$ be the preimage of $a$ in $R$, we
then have

\begin{align*}
	i^{\cT_1}(a_c)& = t^\Phi_{\alpha}(a_c)\\
	& = t^\Phi_{\alpha} \circ i_{0,\alpha}^{\cT_0} (a_R)\\
	& = i^{X_0}(a_R) = i^{Y_0}(a_R)\\
	& = t^\Psi_\beta \circ i_{0,\beta}^{\cU_0} (a_R)\\
	& = t^\Psi_\beta(a_d)\\
	& = i^{\cU_1}(a_d),
\end{align*}
as desired. This proves the claim.
\hfill       $\square$

\medskip
\noindent
{\em Remark.} See \cite{mouse.suslin}[Lemma 2.24] for an argument that overlaps with this one.

\medskip
\noindent
Let us define, for any $\beta < \kappa$,
\[
	f_R^a(\beta) = b \text{ iff } \exists d (R \le d \wedge \beta = \bar{\pi}_{0,\infty}(\kappa_d)
	\wedge b = \bar{\pi}_{0,\infty}(a_d)).
	\]
Claim 2 implies that $f^a_R$ is a function. It is clear that $\dom(f^a_R) \in U$. 

\bigskip
\noindent
{\em Claim 3.} Let $R$ and $S$ be such that $a \in \ran(\pi_{R,\infty})$ and
$a \in \ran(\pi_{S,\infty})$; then $[f^a_R]_U = [f^a_S]_U$.

\medskip
\noindent
{\em Proof.} For $U$ a.e. $\beta$, there is a $d$ such that
$R \le d$, $S \le d$, and $\beta = \kappa_d$. For any such $\beta$
and $d$, $f^a_R(\beta) = \bar{\pi}_{d,\infty}(a_d) = f^a_S(\beta)$.
\hfill    $\square$

\medskip
\noindent
We shall set $\sigma(a) = [f^a_R]_U$. This leads to
$\sigma([a,g]^{M_\infty}_E) = j_U(g)([f^a_R]_U)$, or in
other words, 
\[
	\sigma([a,g]^{M_\infty}_E) = [g \circ f^a_R]_U.
	\]
The following claim implies that this works.

\bigskip
\noindent
{\em Claim 4.} Let $\ult(M_\infty,E) \models
\varphi[[a_0,g_0], ...,[a_n,g_n]]$, and let
$a_i \in \ran(\pi_{R_i,\infty})$ for all $i \le n$.
Then for $U$-a.e. $\beta$, $M_\infty \models
\varphi[g_0(f^{a_0}_{R_0}(\beta)), ..., g_n(f^{a_n}_{R_n}(\beta))].$

\medskip
\noindent
{\em Proof.} Let us assume $n =0$, and write $g = g_0$, $a = a_0$, and $R = R_0$.
By Claim 3, we may assume that $g \in \ran(\pi_{R,\infty})$. 
It is enough to show that whenever
$R \le d$, then $M_\infty \models \varphi[g(\bar{\pi}_{d,\infty}(a_d)].$

Let $g = \pi_{R,\infty}(g_R)$, $a = \pi_{R,\infty}(a_R)$,
and $E = \pi_{R,\infty}(E_R)$. We have that
for $(E_R)_{a_R}$ a.e. $U$, $R \models \varphi[g_R(u)]$.
Now let $R \le d$, and set $g_d = \pi_{R,M_d}(g_R)$.
Again, we have that for $(E_d)_{a_d}$ a.e. $u$,
$M_d \models \varphi[g_d(u)]$. 
It follows that
\[
	\ult(M_d,E_d) \models \varphi[i_{E_d}^{M_d}(g)(a_d)].
\]
Letting $S = \ult(M_d,E_d)$, we have $\pi_{M_d,\infty} =
\pi_{S,\infty}\circ i_{E_d}^{M_d}$, so
\[
M_\infty \models \varphi[g(\pi_{S,\infty}(a_d)].
\]
But $Q_d = S|lh(E_d)$, and $\lambda(E_d)$ is a cutpoint of $S$,
so by strategy coherence
\[
\pi_{S,\infty}\restriction \lambda(E_d) = \bar{\pi}_{d,\infty} \restriction \lambda(E_d).
\]
Thus $M_\infty \models \varphi[g(\bar{\pi}_{d,\infty}(a_d)]$, as desired.
\hfill    $\square$

\bigskip
\noindent
By Claim 4, the map $\sigma([a,g]^{M_\infty}_E) = [g \circ f^a_R]_U$ is well defined
and elementary. Written otherwise, $\sigma(i_E^{M_\infty}(g)(a)) = j_U(g)([f_R^a])$.
Applied to constant functions $g$, this tells us $j_U \restriction M_\infty =
\sigma \circ i_E^{M_\infty}$. Evaluating at $\kappa$, we see that $i_E^{M_\infty}(\kappa)
\le j_U(\kappa)$.

\section{Proof of Theorem \ref{measurebound}}

\begin{definition} $S_\kappa$ is the pointclass of $\kappa$-Suslin sets.
\end{definition}

     Let $\kappa$ be a limit of Suslin cardinals, and $\cf(\kappa)>\omega$. Put
     \[
	     \Delta = \bigcup_{\alpha < \kappa}S_\alpha.
	     \]
$\kappa = \delta(\Delta)$ is the sup of the lengths of prewellorderings in $\Delta$, as well as
its Wadge rank. (See the proof of 3.8 of \cite{jackson}.)
  Let $\Gamma$ be the boldface pointclass such that 
  \[
	  \Delta =  \Gamma \cap \breve{\Gamma} \text{ and $\breve{\Gamma}$ has the Separation property.} 
	    \]
The paper
\cite{closureproperties} shows there is such a $\Gamma$, identifies $\Gamma$ as the class of $\bf{\Sigma}^1_1$-bounded
unions of sets in $\bigcup_{\alpha<\kappa}S_\alpha$, and shows $\forall^\R \Gamma \subseteq \Gamma$. 
Jackson has shown that $\Gamma$ is precisely
the class of all $p[T]$, for $T$ a homogeneous tree on $\omega \times \kappa$,
and that it has the scale property.
See \cite{jackson}[3.8]. It is also shown there that $S_\kappa = \exists^\R \Gamma$.
Another somewhat useful fact is that there is a regular $\Gamma$ norm $\varphi$ on a
complete $\Gamma$ set such that $\le_{\varphi}$ has order type $\kappa$.
(See \cite{jackson}[2.22].) 

Let us fix such a $\Gamma$ norm $\varphi \colon B \to \R$.
Using $\varphi$ and the uniform coding lemma (see \cite{kowo}), we get a coding of
subsets of $\kappa$. To be precise, let $B_\alpha = \lbrace (x,y) \mid \varphi(x) \le \varphi(y) \le \alpha \rbrace$.
For any $A \subset \kappa$, there is a real $x$ and $\Sigma^1_1$ formula $\psi$
such that for all $\alpha < \kappa$ and $y$ such that $\varphi(y) = \alpha$ and $\gamma > \alpha$,
\[
\alpha \in A \Leftrightarrow \psi(y, B_\gamma,x).
\]
($B$ can occur negatively in $\psi$.) We can assume $\psi$ is fixed for all
$x$ by using a universal formula. For any real $x$, let
\[
	\alpha \in A_x \text{ iff } \exists \gamma > \alpha \exists y (\varphi(y) = \alpha
	\wedge \psi(y,B_\gamma,x)).
\]
So $P(\kappa) = \lbrace A_x \mid x \in \R \rbrace$. Using the Godel pairing
we let 
\[
	f_x = \lbrace (\alpha,\beta) \mid (\alpha,\beta) \in A_x \rbrace.
	\]
	Of course, $f_x$ may not be a function. We say $x$ is {\em single valued}
	if $f_x$ is a function. (It need not be total, however.)

	    We define the Martin class, or envelope, of $\Gamma$ by
\begin{align*}
	A \in \Lambda(\Gamma,\kappa) \text{ iff } &  \exists \langle A_\alpha \mid \alpha<\kappa\rangle
			   [ \forall \alpha < \kappa (A_\alpha \in \Delta) \text{ and }\\
			   & \forall^*d \exists \alpha < \kappa (A \cap \lbrace x \mid x \le d \rbrace =
			   A_\alpha \cap \lbrace x \mid x \le d \rbrace).
\end{align*}
The main thing is

\bigskip
\noindent
{\em Claim 1.} Let $U \in \meas_\kappa$, and put 
$x \prec y$ iff ($x$ and $y$ are single valued and defined $U$-a.e.,
and $[f_x]_U \le [f_y]_U$); then 
$\prec$ is in $\Lambda(\Gamma,\kappa)$.

\medskip
\noindent
{\em Proof.} For $\gamma < \beta < \kappa$, put
\begin{align*}
	(x,y) \in A_{\beta,\gamma} \Leftrightarrow & (\gamma \in \dom(f_x) \cap \dom(f_y) \wedge
	f_x(\gamma) \le f_y(\gamma) <\beta \\
	   & \wedge f_x \cap \beta \times \beta \text{ and } f_y \cap \beta \times \beta \text{ are single valued.}
\end{align*}
It is enough to show that for any $d$,  there is a $(\beta,\gamma)$ such that
$\prec$ agrees with $A_{\beta,\gamma}$ on the reals $\le d$. But by countable completeness,
we can find $\gamma$ such that for all single valued $x \le d$,
$\gamma \in \dom(f_x)$ iff $\dom(f_x) \in U$. Similarly, we can
arrange that for $x,y \le d$ single-valued with domains in $U$, $f_x(\gamma) \le f_y(\gamma)$
iff $[f_x]_U \le [f_y]_U$. Finally, since $\kappa$ has uncountable cofinality,
we can choose $\beta$ large enough that all relevant $f_x(\gamma)$ are below
$\beta$, and any non-single-valued $x \le \beta$ are such that
$f_x \cap \beta \times \beta$ is not single valued. This proves the Claim.
\hfill    $\square$

\medskip

     Let $\lambda$ be the least Suslin cardinal $>\kappa$. As shown in
\cite{jackson}, the universal $\breve{S}_\kappa$ set has a semi-scale
all of whose norm relations are in the envelope $\Lambda(\Gamma,\kappa)$.\footnote{ Jackson
and Woodin showed there is a self-justifying system sealing the envelope, in fact.}

If $S_\kappa$ is closed under $\forall^{\R}$, or equivalently $S_\kappa = \Gamma$,
we get $\Lambda(\Gamma,\kappa)$ is closed under real quantifiers. Martin's
non-uniformizability result then shows that
$\lambda$ is at least prewellordering ordinal of $\Lambda(\Gamma,\kappa)$.
( See \cite{jackson}[3.17].) Combined
with Claim 1, this gives $\lambda \ge \sup(\lbrace j_U(\kappa) \mid
U \in \meas_\kappa \rbrace)$. So $\lambda$, this sup, and the prewellordering ordinal
of $\Lambda$ coincide.

   Now let us assume that $\Gamma$ is not closed under $\exists^{\R}$, and look at the
    projective-like hierarchy above it.
We write $\Pi_1 = \Gamma$, $\Sigma_1 = \breve{\Gamma}$, and $\Pi_{n+1} = \forall^\R \Sigma_n$
and $\Sigma_{n+1} = \exists^\R \Pi_n$. For $n>1$, these pointclasses have the usual
closure properties of the levels of the projective hierarchy. By periodicity,
the $\Pi_{2n+1}$ and $\Sigma_{2n+2}$ have the scale property.
Since we are asuming $S_\kappa \neq \Pi_1$, we get $S_\kappa = \Sigma_2$. It follows
by the ordinary projective hierarchy arguments that $\lambda$ has cofinality
$\omega$, $S_\lambda = \Sigma_3$, and $\lambda^+$ is the next Suslin after
$\lambda$, and the prewellordering ordinal of $\Pi_3$.

\medskip
\noindent
{\em Remark.} In the present case, $S_\kappa = \Sigma_2$ is the class of all
$\kappa$-length unions of sets in $\Delta$. It is therefore properly 
included in $\Lambda(\Gamma,\kappa)$. What Martin's proof shows is
that there is a $\Pi_2$ relation with no uniformization in
$\Lambda(\Gamma,\kappa)$.

\medskip

     But now let $\prec$ be any prewellorder in $\Lambda(\Gamma,\kappa)$.
If $\lambda^+ \le |\prec|$, then $\prec$ is not $\lambda$-Suslin
by Kunen-Martin. It follows by Wadge that a universal $\Pi_3$
set is Wadge below $\prec$.  But we can uniformize every $\Pi_2$
relation in $\Pi_3$, since the latter has the scale property. This
contradicts Martin's theorem. We have proved that
$|j_U(\kappa)| \le \lambda$ for all $U \in \meas_\kappa$.
(But the sup might be $\lambda^+$.)

      Now let us assume $\kappa$ is regular. This makes $\Gamma$ a nicer pointclass,
      closed under countable unions and intersections, for example.
Also, $\kappa$ has the strong partition property from the usual arguments using the uniform
coding lemma. A general fact is that if $\kappa$ is any cardinal with the strong partition
if $U$ is semi-normal, that is gives every club set measure one, then $j_U(\kappa)$ is regular.
So, if $j_U(\kappa)
\geq \lambda$, then $j_U(\kappa)\geq \lambda^+$. We have just shown this is not the case.

      Finally, suppose that $\kappa$ is singular. Let $W$ be the $\omega$-club ultrafilter
      on $\cf(\kappa)$, which exists because we are in the range of the $\hd$-analysis.
      Jackson shows that $j_{U \times W}(\kappa)$ is a cardinal. ( Proof to come!)
      This shows $j_{U \times W}(\kappa) < \lambda$, so $j_U(\kappa)<\lambda$.

This completes the proof of Theorem \ref{measurebound}.

\section{ Proofs of \ref{twoimpliesone}, \ref{addsauce}, and \ref{cofomegalimit}}

\bigskip
\noindent
{\em Proof of \ref{twoimpliesone}.}

    We are given a projectum stable mouse pair $(P,\Sigma)$, and $\kappa < o(M_\infty(P,\Sigma)$
such that $\kappa$ is a cardinal of $V$, and a cutpoint of $M_\infty(P,\Sigma)$.
By \cite{mouse.suslin}[2.19], $|\tau_\infty(P,\Sigma)|$ is a Suslin cardinal,
and $|\tau_\infty(P,\Sigma)| = |o(M_\infty(P,\Sigma)| \ge \kappa$. So if
$\tau_\infty(P,\Sigma) < (\kappa^+)^V$, then $\kappa$ is a Suslin cardinal,
as desired. So assume $(\kappa^+)^V \le \tau_\infty(P,\Sigma)$.

   Set $M_\infty = M_\infty(P,\Sigma)$. Suppose that $o(\kappa)^{M_\infty} \ge 
   \tau_{\infty}(P,\Sigma)$.
Since $\kappa$ is a cutpoint of $M_\infty$, this implies
that  $(P,\Sigma)$ has a top block, and $\beta_\infty(P,\Sigma) = \kappa$.
Then by \cite{mouse.suslin}[2.27], $\kappa$ is a Suslin cardinal, as desired.

   So suppose that $o(\kappa)^{M_\infty} < \tau_\infty(P,\Sigma)$.
The following little lemma is useful.

\begin{lemma}\label{cutpointlemma} Let $(R,\Omega)$ be a projectum stable mouse pair,
	$R_\infty = M_\infty(R,\Omega)$, and
	$k = k(R)$. Let 
	\[
	\gamma =\sup\lbrace \eta,o(\eta)^{R_\infty}\rbrace,
	\]
where $\eta$ is a cardinal of $R_\infty$, and
	\[
		\xi = (\gamma^+)^{R_\infty}.
		\]
Suppose $\xi \le \rho_k(R_\infty)$;
then there is a $(Q,\Psi)$ such that
	\begin{itemize}
		\item[(a)]
	$M_\infty(Q,\Psi) = R_\infty|\xi$, and
	\item[(b)] 
$\tau_\infty(Q,\Psi) \le \gamma$.
	\end{itemize}
Thus $|\gamma|$ is a Suslin cardinal.
\end{lemma}
\begin{proof}
By replacing $(R,\Omega)$ with an iterate of itself, we may assume
that we have $\bar{\eta},\bar{\gamma},$ and $\bar{\xi}$
such that $\pi_{R,\infty}(\langle \bar{\eta},\bar{\gamma},\bar{\xi} \rangle)
= \langle \eta,\gamma,\xi \rangle$. By coherence, $\bar{\xi}$
is a cutpoint of $R$. Also, $\bar{\xi} \le \rho_k(R)$. $\bar{\xi}$
is regular in $R$, so if it were not $R$-regular, we would have
$\bar{\xi} = \rho_k(R)$ and some $\Sigma_k(R)$ partial $f$
with $\dom(f) \subseteq \bar{\gamma}$ and $\ran(f)$ cofinal
in $\bar{\xi}$. This easily yields $\rho_k(R) \le \bar{\gamma}$,
contradiction. Thus $\bar{\xi}$ is $R$-regular.

But then we can take $Q = R|\bar{\xi}$ and $\Psi = \Omega_Q$.
\end{proof}

	\medskip
	\noindent
	Now let $\gamma = o(\kappa)^{M_\infty}$, and $\xi = (\gamma^+)^{M_\infty}$.
	We are assuming $\gamma<\tau_\infty(P,\Sigma)$, and
	$\tau_\infty(P,\Sigma) < \rho_k(M_\infty)$ by its definition,
	so $\xi \le \rho_k(M_\infty)$. Thus by the lemma,
	$|\gamma|$ is a Suslin cardinal, so we may assume $(\kappa^+)^V \le \gamma$,
	otherwise we're done. Let then $(Q,\Psi)$ be such that
	$M_\infty(Q,\Psi) = M_\infty|\xi$. It is easy to see that
	$\kappa = \beta_\infty(Q,\Psi)$. Thus by \cite{mouse.suslin}[2.27],
	$\kappa$ is a Suslin cardinal. This completes the proof
	of \ref{twoimpliesone}.
	\hfill   $\square$

	\bigskip
	\noindent
	We turn to \ref{addsauce} and \ref{cofomegalimit}.

    Let $(P,\Sigma)$ be a projectum stable mouse pair, and $\kappa$ a limit of Suslin cardinals,
    and $\kappa < o(M_\infty)$ where $M_\infty = M_\infty(P,\Sigma)$. 
Replacing $(P,\Sigma)$ with
an iterate of itself, we may assume $\kappa = \pi_{P,\infty}(\kappa_P$.
Let $\pi_{P,\infty}(\kappa_P)=\kappa$.
 
    That $\kappa$ is a limit of cutpoints in $P$ was proved in \cite{mouse.suslin}.
(Let $\mu$ be least such that $\mu < \kappa$ and $o(\mu)^{M_\infty} \ge \kappa$.
By Cor. 2.42 of \cite{mouse.suslin}, there are no Suslin cardinals strictly between
$\mu$ and $o(\mu)^{M_\infty}$. Contradiction.) We can also prove it using the
measure existence result \ref{paulgrigor},
and a softer coarser form of the measure bounding result \ref{measurebound}.

\bigskip
\noindent
{\em Proof of \ref{addsauce}.}

\medskip
\noindent
	We assume $\cf(\kappa)>\omega$ and $\kappa^+ \le o(M_\infty)$.
Thus $\kappa^+ \le  \tau_\infty(P,\Sigma)$.  Let $\lambda$ be the least Suslin cardinal $>\kappa$,
so that by $\cite{jackson}[\S3]$, $\cf(\lambda)=\omega$, and hence
$\kappa^+ < \lambda$. Since $|\tau_\infty(P,\Sigma)|$ is a Suslin
cardinal, 
$ \lambda \le \tau_\infty(P,\Sigma)$. Thus $\lambda \le \rho_k(M_\infty)$,
where $k = k(P)$.

  By Lemma \ref{cutpointlemma}, there is no cutpoint
$\xi$ of $M_\infty$ such that $\xi = (\gamma^+)^{M_\infty}$
for some $\gamma$, and $\kappa^+ < \xi <\lambda$. For otherwise,
there would be Suslin cardinals in the interval $(\kappa,\lambda)$.
It follows that there is a $\mu < \kappa^+$ such that
$o(\mu)^{M_\infty} \ge \lambda$. Since $\kappa$ is a cutpoint,
$\kappa \le \mu$.

 By coherence, we get that if $(\kappa^+)^V \le o(\kappa)^{M_\infty}$,
 then $\lambda \le o(\kappa)^{M_\infty}$. 

 Let $\mu < \kappa^+$ be such that $o(\mu)^{M_\infty} \ge \lambda$.
Let $E$ be a total $M_\infty$ extender with critical point $\mu$.
By Theorem \ref{paulgrigor}, there is an ultrafilter $U$ on $\mu$
such that $\lambda(E) \le j_U(\mu)$. But $|\mu| = \kappa$, so by
Theorem \ref{measurebound}, $j_U(\kappa) < \lambda$, and thus
$j_U(\mu) < \lambda$. So $o(\mu)^{M_\infty} \le \lambda$, hence
$o(\mu)^{M_\infty} = \lambda$.

Finally, suppose $S_\kappa$ is closed under $\forall^\R$. We must see
$o(\kappa)^{M_\infty} \ge \kappa^+$. If not, by Lemma
\ref{cutpointlemma} we get $(Q,\Psi)$ such that
$\kappa < o(M_\infty(Q,\Psi)) < \kappa^+$. Thus $\cd(\Psi)$
is $\kappa$-Suslin. But $\Psi$ is a complete strategy, so
\[
	\cT \text{ is not by $\Psi$ iff } \exists \alpha < \lh(\cT) (\Psi(\cT \restriction \alpha)
	\neq [0,\alpha)_T).
\]
Thus $\neg \cd(\Psi)$ is also $\kappa$-Suslin. Since $S_\kappa$ is inductive
like, we get that $\cd(\Psi) \in S_\alpha$ for some $\alpha < \kappa$,
contrary to Kunen-Martin and the fact that $\kappa \le o(M_\infty(Q,\Psi)$.
\hfill  $\square$

\bigskip
\noindent
{\em Proof of \ref{cofomegalimit}}

\medskip
\noindent
$\kappa < \rho_k(M_\infty)$ because we demanded $(\kappa^+)^V \le
o(M_\infty)$. So if $\kappa$ were measurable in $M_\infty$, the
set of images of iteration points of the order zero measure would have
uncountable cofinality, contrary to $\cf(\kappa)=\omega$. $\kappa^+ =
(\kappa^+)^V$ is regular, in fact measurable, in $V$ because it is the
prewellordering ordinal of a $\Pi^1_1$-like pointclass. So
$\kappa^+ = o(M_\infty)$ is impossible, as $\pi_{P,\infty}``o(P)$
would be cofinal in $\kappa^+$. A similar argument shows $\kappa^+$ is
measurable in $M_\infty$.

   Finally, suppose toward contradiction that $\kappa < \mu < \kappa^+$,
   and $\mu$ is a cutpoint of $M_\infty$ such that $o(\mu)^{M_\infty} \ge \kappa^+$.
Then $\cf(\mu) > \omega$ because $\mu$ is measurable in $M_\infty$ by a total measure,
and $o(\mu)^{M_\infty} > \kappa^+$ by coherence and the fact that $\kappa^+$ is
measurable in $M_\infty$..
Applying the result of Sargsyan, out in $V$ there is an ultrafilter
$U$ on $\mu$ such that $j_U(\mu) > \kappa^+$.  But $|\mu| = \kappa$, so we
can use a bijection to replace $U$ with an ultrafilter $W$ on $\kappa$
such that $j_W = j_U$. Since $\cf(\kappa)=\omega$, we may assme that $W$
is actually an ultrafilter on some $\eta < \kappa$.

    $\kappa$ is a limit of Suslin cardinals, so easy measure bounding
    gives $j_W(\kappa) = \kappa$. If $f$ maps $\kappa$ onto $\mu$,
    then $j_W(f)$ maps $\kappa$ onto $j_W(\mu)$, so $j_W(\mu)<\kappa^+$,
    contradiction.
    \hfill  $\square$

    \bigskip

    We omit the proof of \ref{conjectureprojective}. It is like the proofs above,
    but uses Jackson's measure-bounding results for measures on projective ordnals.

    \end{document}